\documentclass[11pt]{amsart}
\usepackage{amsmath,amsthm,amsfonts,amssymb,latexsym,epsfig}

\newtheorem{theorem}{Theorem}[section]

\newtheorem{lemma}{Lemma}[section]

\newtheorem{cor}{Corollary}[section]

\numberwithin{equation}{section}

\theoremstyle{definition}

\theoremstyle{remark}

\begin{document}
\title{On $l^p$ norms of weighted mean matrices}
\author{Peng Gao}
\address{Division of Mathematical Sciences, School of Physical and Mathematical Sciences,
Nanyang Technological University, 637371 Singapore}
\email{penggao@ntu.edu.sg}
\date{August 31, 2008.}
\subjclass[2000]{Primary 47A30} \keywords{Carleman's inequality, Hardy's inequality, weighted mean matrices}


\begin{abstract}
 We study $l^{p}$ operator norms of weighted mean matrices using the approaches of Kaluza-Szeg\"o and Redheffer.
 As an application, we prove a conjecture of Bennett.
\end{abstract}

\maketitle
\section{Introduction}
\label{sec 1} \setcounter{equation}{0}

  Suppose throughout that $p\neq 0, \frac{1}{p}+\frac{1}{q}=1$.
  For $p \geq 1$, let $l^p$ be the Banach space of all complex sequences ${\bf a}=(a_n)_{n \geq 1}$ with norm
\begin{equation*}
   ||{\bf a}||_p: =(\sum_{n=1}^{\infty}|a_n|^p)^{1/p} < \infty.
\end{equation*}
  The celebrated
   Hardy's inequality (\cite[Theorem 326]{HLP}) asserts that for $p>1$,
\begin{equation}
\label{eq:1} \sum^{\infty}_{n=1}\Big{|}\frac {1}{n}
\sum^n_{k=1}a_k\Big{|}^p \leq \Big (\frac
{p}{p-1} \Big )^p\sum^\infty_{n=1}|a_n|^p.
\end{equation}
   Hardy's inequality can be regarded as a special case of the
   following inequality:
\begin{equation}
\label{01}
   \Big | \Big |C \cdot {\bf a}\Big | \Big |^p_p =\sum^{\infty}_{n=1} \Big{|}\sum^{\infty}_{k=1}c_{n,k}a_k
    \Big{|}^p \leq U \sum^{\infty}_{n=1}|a_n|^p,
\end{equation}
   in which $C=(c_{n,k})$ and the parameter $p>1$ are assumed
   fixed, and the estimate is to hold for all complex
   sequences ${\bf a} \in l^p$. The $l^{p}$ operator norm of $C$ is
   then defined as
\begin{equation*}
\label{02}
    ||C||_{p,p}=\sup_{||{\bf a}||_p \leq 1}\Big | \Big |C \cdot {\bf a}\Big | \Big |_p.
\end{equation*}
   It follows that inequality \eqref{01} holds for any ${\bf a} \in l^p$ when $U^{1/p} \geq ||C||_{p,p}$ and fails to hold for some ${\bf a} \in l^p$
   when $U^{1/p} <||C||_{p,p}$.
    Hardy's inequality thus asserts that the Ces\'aro matrix
    operator $C$, given by $c_{n,k}=1/n , k\leq n$ and $0$
    otherwise, is bounded on {\it $l^p$} and has norm $\leq
    p/(p-1)$. (The norm is in fact $p/(p-1)$.)

    We say a matrix $A=(a_{n,k})$ is a lower triangular matrix if $a_{n,k}=0$ for $n<k$ and a lower triangular matrix $A$ is a summability matrix if
    $a_{n,k} \geq 0$ and
    $\sum^n_{k=1}a_{n,k}=1$. We say a summability matrix $A$ is a weighted
    mean matrix if its entries satisfy:
\begin{equation}
\label{021}
    a_{n,k}=\lambda_k/\Lambda_n,  ~~ 1 \leq k \leq
    n; \hspace{0.1in} \Lambda_n=\sum^n_{i=1}\lambda_i, \lambda_i \geq 0, \lambda_1>0.
\end{equation}

    Hardy's inequality \eqref{eq:1} now motivates one to
    determine the $l^{p}$ operator norm of an arbitrary summability matrix $A$.
   In an unpublished dissertation \cite{Car}, Cartlidge studied
weighted mean matrices as operators on $l^p$ and obtained the
following result (see also \cite[p. 416, Theorem C]{B1}):
\begin{theorem}
\label{thm02}
    Let $1<p<\infty$ be fixed. Let $A$ be a weighted mean matrix given by
    \eqref{021}. If
\begin{equation}
\label{022}
    L=\sup_n\Big(\frac {\Lambda_{n+1}}{\lambda_{n+1}}-\frac
    {\Lambda_n}{\lambda_n}\Big) < p ~~,
\end{equation}
    then
    $||A||_{p,p} \leq p/(p-L)$.
\end{theorem}

   There are several published proofs of Cartlidge's
   result. Borwein \cite{Bor} proved a far more general
result than Theorem \ref{thm02} on the $l^p$ norms of generalized
Hausdorff matrices. Rhoades \cite[Theorem 1]{Rh} obtained a slightly
general result
 than Theorem \ref{thm02}, using a modification of the proof of
 Cartlidge. Recently, the author \cite{G1} also gave a simple proof of Theorem
 \ref{thm02}.

   We note here that by a change of variables $a_k \rightarrow a^{1/p}_k$ in \eqref{eq:1} and on letting $p \rightarrow +\infty$, one obtains the
   following well-known Carleman's inequality  \cite{Carlman}, which asserts that for convergent infinite series $\sum a_n$ with non-negative terms,
   one has
\begin{equation*}
   \sum^\infty_{n=1}(\prod^n_{k=1}a_k)^{\frac 1{n}}
\leq e\sum^\infty_{n=1}a_n,
\end{equation*}
   with the constant $e$ being best possible.

   It is then natural to study the following weighted version of Carleman's inequality:
\begin{equation}
\label{1}
   \sum^\infty_{n=1}\Big( \prod^n_{k=1}a^{\lambda_k/\Lambda_n}_k \Big )
\leq E\sum^\infty_{n=1}a_n,
\end{equation}
  where the notations are as in \eqref{021}. The task here is to determine the best constant $E$ so that inequality \eqref{1} holds
  for any convergent infinite series $\sum a_n$ with non-negative terms.
  Note that Cartlidge's result (Theorem \ref{thm02}) implies that when \eqref{022} is satisfied, then for any ${\bf a} \in l^p$, one has
\begin{equation}
\label{4}
   \sum^{\infty}_{n=1}\Big{|}\sum^{n}_{k=1}\frac {\lambda_ka_k}{\Lambda_n}
   \Big{|}^p \leq \Big ( \frac {p}{p-L} \Big )^p \sum^{\infty}_{n=1}|a_n|^p.
\end{equation}
  Similar to our discussions above, by a change of variables $a_k \rightarrow a^{1/p}_k$ in \eqref{4} and on letting $p \rightarrow +\infty$,
  one obtains inequality \eqref{1} with $E=e^{L}$ as long as \eqref{022} is satisfied with $p$ replaced by $+\infty$ there.

  It is our goal in this paper to extend the result of Theorem \ref{thm02} and to study weighted Carleman's inequality as well.
  Our treatment of this subject will not be satisfactory if we don't mention the integral cases of Hardy's inequality as the later often supply
  the motivations for considerations of various discrete cases.
  For this reason, we will first give a brief discussion in Section \ref{sec 1'} on
  the integral Hardy-type inequalities and point out the ties between them and the discrete cases.
  For the rest part of the paper, we will focus ourselves on the $l^p$ spaces as well as their variations.
  A general method towards establishing $l^{p}$ operator norms of weighted mean matrices would be via Carleman's approach,
  which is essentially a use of Lagrange multipliers as we shall explain in details in Section \ref{sec 2}.
  However, this approach is more technically involved so we are looking for other methods that can be used to achieve our goal
  in this paper while technically simpler compared to Carleman's approach.
  Among the many different proofs of Hardy's inequality \eqref{eq:1} as well as its generalizations and extensions in the literature,
  there are notably Kaluza and Szeg\"o's approach \cite{K&S} (see also \cite{K}) and Redheffer's ``recurrent inequalities" \cite{R1}.
It is shown in \cite{G3} that these two methods above are
essentially the same (in \cite{G3}, we credited the approach of
Kaluza and Szeg\"o to Knopp but apparently the paper \cite{K&S} is
earlier) and we shall further show in this paper that Kaluza and
Szeg\"o's approach can be regarded as an approximation to
Carleman's approach in Section \ref{sec 2}. Hence instead of
Carleman's approach, there is not much lost using Kaluza and
Szeg\"o's or Redheffer's approach when studying Hardy-type
inequalities, yet technically they are  much easier to handle.

   In this paper, we shall use Kaluza and Szeg\"o's approach to prove the following extension of Theorem \ref{thm02} (we note here the case $n=1$
of \eqref{022} implies $L > 0$) in Section \ref{sec 3}:
\begin{theorem}
\label{thm03}
    Let $1<p<\infty$ be fixed. Let $A$ be a weighted mean matrix given by
    \eqref{021}. If for any integer $n \geq 1$, there exists a positive constant
    $0<L<p$ such that
\begin{equation}
\label{024}
    \frac {\Lambda_{n+1}}{\lambda_{n+1}} \leq \frac
    {\Lambda_n}{\lambda_n}  \Big (1- \frac
    {L\lambda_n}{p\Lambda_n} \Big )^{1-p}+\frac {L}{p}~~,
\end{equation}
    then
    $||A||_{p,p} \leq p/(p-L)$.
\end{theorem}

    Note that for $p>1$,
 \begin{equation*}
     \Big (1-\frac
    {L\lambda_n}{p\Lambda_n} \Big )^{1-p}  \geq 1+(1-\frac {1}{p})\frac
    {L\lambda_n}{\Lambda_n} + (1-\frac {1}{p}) \frac
    {\lambda^2_n}{\Lambda^2_n} \frac {L^2}{2}.
 \end{equation*}
     It follows from this and \eqref{024} that we have the following
\begin{cor}
\label{cor03}
    Let $1<p<\infty$ be fixed. Let $A$ be a weighted mean matrix given by
    \eqref{021}. If for any integer $n \geq 1$, there exists a positive constant
    $0<L<p$ such that
\begin{equation*}
    \frac {\Lambda_{n+1}}{\lambda_{n+1}}-\frac
    {\Lambda_n}{\lambda_n} \leq L+ \Big ( \frac
    {\lambda_n}{2\Lambda_n} \Big )\Big(1-\frac 1{p}\Big)L^2~~,
\end{equation*}
    then
    $||A||_{p,p} \leq p/(p-L)$.
\end{cor}

   An interesting proof of Hardy's inequality \eqref{eq:1} for $p=2$ is given by Wang and Yuan in \cite{W&Y}.
   Their method interprets the left-hand side of \eqref{eq:1} when $p=2$ as a quadratic form
   so that Hardy's inequality follows from estimations of the
   eigenvalues of the corresponding matrix associated to the quadratic form.
   We point out here that this approach, as we will show in Section \ref{sec 2'},
   can be viewed as an approach via the duality principle of linear operators.
   We will then use the method of Wang and Yuan to give another proof of Theorem \ref{thm03} in Section \ref{sec 2'}.
   As we shall also see there, this approach also has ties with Carleman's approach.

   We point out here Theorem \ref{thm02} can be regarded as the case $p \rightarrow 1^+$ of Theorem \ref{thm03}
   while the case $p \rightarrow +\infty$ of Theorem \ref{thm03} suggests the following result:
\begin{theorem}
\label{thm1}
  Suppose that
\begin{equation}
\label{5}
  M=\sup_n\frac
    {\Lambda_n}{\lambda_n}\log \Big(\frac {\Lambda_{n+1}/\lambda_{n+1}}{\Lambda_n/\lambda_n} \Big ) < +\infty,
\end{equation}
  then inequality \eqref{1} holds with $E=e^M$.
\end{theorem}
   It is easy to see that $M \leq L$ for $L$ defined by \eqref{022} and hence the above theorem does provide a better result
   than what one can infer from Cartlidge's result as discussed above if \eqref{5} were true.
   In fact, a even stronger result exists, namely, the following nice result of Bennett \cite{Be1} (see the proof of Theorem 13):
\begin{theorem}
\label{thm6.1}
  Inequality \eqref{1} holds with
\begin{equation*}
  E=\sup_n\frac
    {\Lambda_{n+1}}{\lambda_{n+1}} \prod^n_{k=1} \Big(\frac {\lambda_{k}}{\Lambda_{k}} \Big)^{\lambda_k/\Lambda_n}.
\end{equation*}
\end{theorem}
   It's shown in \cite{Be1} that the above theorem implies Theorem
   \ref{thm1} and the author has recently obtained the $l^p$ analogues of
   the above theorem which in fact implies Theorem \ref{thm03}. We
   shall not go further here and refer the interested
   reader to the paper \cite{G7} for the details.

   To give applications of Theorem \ref{thm03} or Corollary \ref{cor03}, we note that the following two
     inequalities were claimed to hold for any ${\bf a} \in l^p$ by Bennett ( \cite[p. 40-41]{B4}; see also \cite[p. 407]{B5}):
\begin{eqnarray}
\label{7}
   \sum^{\infty}_{n=1}\Big{|}\frac
1{n^{\alpha}}\sum^n_{i=1}(i^{\alpha}-(i-1)^{\alpha})a_i\Big{|}^p &
\leq & \Big( \frac {\alpha p}{\alpha p-1} \Big )^p\sum^{\infty}_{n=1}|a_n|^p, \\
\label{8}
   \sum^{\infty}_{n=1}\Big{|}\frac
1{\sum^n_{i=1}i^{\alpha-1}}\sum^n_{i=1}i^{\alpha-1}a_i\Big{|}^p &
\leq & \Big(\frac {\alpha p}{\alpha p-1} \Big
)^p\sum^{\infty}_{n=1}|a_n|^p,
\end{eqnarray}
     whenever $p>1, \alpha p >1$. We note here the constant $(\alpha p /(\alpha p-1))^p$ is best possible (see \cite{Be1}).

No proofs of the above two inequalities were supplied in \cite{B4}-\cite{B5}. The author \cite{G} and Bennett \cite{Be1}
     proved inequalities \eqref{7} for $p>1, \alpha \geq 1, \alpha p >1$ and
     \eqref{8} for $p>1, \alpha \geq 2$ or $0< \alpha \leq 1, \alpha p >1$
     independently. The proofs of \eqref{8} in both \cite{G} and \cite{Be1} are the same, they both use the result of Cartlidge (Theorem \ref{thm02} above).
Recently, the author \cite{G3} has shown that inequalities \eqref{8} hold for $p \geq 2, 1 \leq
   \alpha \leq 1+1/p$ or $1 < p \leq 4/3, 1+1/p \leq
   \alpha \leq 2$.

   In connection to \eqref{8}, Bennett \cite[p. 829]{Be1} further conjectured that inequality \eqref{1} holds
   for $\lambda_k=k^{\alpha}$ for $\alpha > -1$ with $E=e^{1/(\alpha+1)}$. As the cases $-1 < \alpha \leq 0$ or $\alpha \geq 1$
   follow directly from the know cases of inequalities \eqref{8} upon changes of
   variables $\alpha \rightarrow \alpha+1, a_k \rightarrow a^{1/p}_k$ and on letting $p \rightarrow +\infty$, the only cases left unknown are when $0< \alpha <1$. As an application of Theorem \ref{thm03}, we shall show that \eqref{8} hold for $ p \geq 2$ and $1 \leq \alpha \leq 2$ in Section \ref{sec 4} and this will in turn imply that Bennett's conjecture is true by our discussions above.

\section{Integral Hardy-type inequalities and Their Discrete Analogues}
\label{sec 1'} \setcounter{equation}{0}

  When $p>1$, the (one dimension) integral Hardy-type inequalities are integral inequalities of the form
\begin{equation}
\label{0}
  \int^b_a s(x)|y(x)|^pdx \leq  \int^b_ar(x)|y'(x)|^pdx,
\end{equation}
  where $r,s$ are non-negative measurable functions on $(a,b)$ and $y(x)$ absolutely continuous on $(a,b)$ subject to certain boundary conditions ($y(a)=0$ or $y(b)=0$ or both) with $r(x)|y'(x)|^p$ integrable on $(a,b)$. The classical integral Hardy's inequality \cite[Theorem 327]{HLP} corresponds to the case of \eqref{0} with $a=0, b=+\infty$, $r=q^p, s(x)=x^{-p}$, $y(x)=\int^x_0f(t)dt$ with $f(x) \geq 0$ and $f(x) \in L^p(0, +\infty)$. Note that in this case the boundary condition $y(0)=0$ is satisfied.

  One can certainly consider analogues of \eqref{0} for other $p$'s. We point out here the above mentioned classical integral Hardy's inequality still holds when one replaces $p>1$ by $p<0$ and this is a result of Beesack \cite[(3.2.13)]{Bee1}. In the case of $0<p<1$, the classical integral Hardy's inequality holds \cite[Theorem 337]{HLP} with inequality reversed when one replaces $r(x)=q^p$ by $r(x)=(-q)^p$ and $y(x)=\int^x_0f(t)dt$ by $y(x)=\int^{\infty}_{x}f(t)dt$ (so the $y(\infty)=0$ in this case).

   By a change of variables $x \rightarrow x^{1+\alpha q}, y(x^{1+\alpha q}) \rightarrow y(x)$ when $p>1$ and $1+\alpha q>0$, we can rewrite the classical integral Hardy's inequality as
\begin{equation*}
  \int^{\infty}_0 \frac {|y(x)|^p}{x^{(1+\alpha)p}}dx \leq \Big (\frac {q}{1+\alpha q} \Big)^p\int^{\infty}_0 \frac {|y'(x)|^p}{x^{\alpha p}}dx.
\end{equation*}
  We point out here and in what follows in this section, unless otherwise mentioned, all the constants are best possible.

  More generally, one has the following result for $p \neq 0$ (with inequalities reversed when $0<p<1$):
\begin{equation}
\label{1.1}
   \int^{\infty}_0 \frac {|y(x)|^p}{x^{(1+\alpha)p}}dx \leq \Big |\frac {q}{1+\alpha q} \Big|^p\int^{\infty}_0 \frac {|y'(x)|^p}{x^{\alpha p}}dx.
\end{equation}
  The $p>1$ case above is \cite[Theorem 330]{HLP} and here we require $y(0)=0$ when $(1+\alpha)p>1$ and $y(\infty)=0$ when $(1+\alpha) p<1$. The $p<0$ case above is \cite[(3.2.13)]{Bee1} and here we require $y(\infty)=0$ when $(1+\alpha)p>1$ and $y(0)=0$ when $(1+\alpha) p<1$. The case $0<p<1$ above is \cite[Theorem 347]{HLP} and here we require $y(0)=0$ when $(1+\alpha)p>1$ and $y(\infty)=0$ when $(1+\alpha) p<1$. If we now write $y(x)=\int^x_0f(t)t^{\alpha}dt$ when $y(0)=0$ or $y(x)=\int^{\infty}_{x}f(t)t^{\alpha}dt$ when $y(\infty)=0$ with $f(x) >0$, then we can recast inequalities \eqref{1.1} as
\begin{equation}
\label{1.2}
    \int^{\infty}_{0}\Big ( \frac {\int^x_0f(t)t^{\alpha}dt}{\int^x_0t^{\alpha}dt} \Big )^pdx \leq \Big |\frac {(1+\alpha )q}{1+\alpha q} \Big|^p\int^{\infty}_{0}f^p(x)dx,
\end{equation}
   for the cases $p>1, (1+\alpha)p>1$ or $p<0, 1+\alpha>0$. The reversed inequalities also hold for $0<p<1, (1+\alpha)p>1$. Similarly, we have
\begin{equation}
\label{1.4}
    \int^{\infty}_{0}\frac {\Big ( \int^{\infty}_xf(t)t^{\alpha}dt \Big )^p}{x^{(1+\alpha)p}}dx \leq \Big |\frac {q}{1+\alpha q} \Big|^p\int^{\infty}_{0}f^p(x)dx,
\end{equation}
   for the case $p>1, 0<1+\alpha<1/p$ with the above inequalities reversed for $0<p<1, 0<1+\alpha<1/p$ and
\begin{equation*}
   \int^{\infty}_{0}\frac {\Big ( \int^{x}_0f(t)t^{\alpha}dt \Big )^p}{x^{(1+\alpha)p}}dx \leq \Big |\frac {q}{1+\alpha q} \Big|^p\int^{\infty}_{0}f^p(x)dx,
\end{equation*}
  for the case $p<0, 1/p<1+\alpha<0$ and
\begin{equation*}
    \int^{\infty}_{0}\Big ( \frac {\int^{\infty}_xf(t)t^{\alpha}dt}{\int^{\infty}_xt^{\alpha}dt} \Big )^pdx \leq \Big |\frac {(1+\alpha)q}{1+\alpha q} \Big|^p\int^{\infty}_{0}f^p(x)dx,
\end{equation*}
   for the cases $p>1, 1+\alpha<0$ or $p<0, 1+\alpha<1/p$ with the above inequalities reversed for $0<p<1, 1+\alpha<0$.

  There have been extensive studies on the integral Hardy-type inequalities in the literature. We refer the reader to the articles \cite{H&L1}, \cite{Bee1}, \cite{Boyd}, \cite{Bee2}, \cite{S} and the references therein for more details.
  One possible way of generalizing the above inequalities is to consider integrations with respect to other measures other than the Lebesgue measure. For example, using the counting measure for the case $\alpha=0$ of \eqref{1.2} leads to the discrete Hardy's inequality \eqref{eq:1} for $p>1$. From this point of view, we see that inequalities \eqref{8} correspond to the general cases of \eqref{1.2}.

   As another example, we note that Levin and Ste\v ckin
  \cite[Theorem D.61]{L&S} proved that for $0 < p  \leq 1/3$,
\begin{equation}
\label{1.6'}
  \sum^{\infty}_{n=1}\Big( \frac 1{n} \sum^{\infty}_{k=n}a_k \Big
  )^p \geq \Big(\frac {p}{1-p} \Big
)^p \sum^{\infty}_{n=1}a^p_n,
\end{equation}
  where and from now on in this section, we assume $a_n \geq 0$ for all $n$.
  Their result can be viewed as the discrete analogues to the case $p<0, \alpha=1$ of \eqref{1.4} and this result improved upon that of Theorem 345 of \cite{HLP}, which asserts the above inequalities hold with smaller constants $p^p$ for $0<p<1$.
Recently, the author has given a simple proof \cite{G3} of the result
of Levin and Ste\v ckin and also extended their result to a case where $p$ is slightly bigger than $1/3$.

   In general, integral Hardy-type inequalities suggest that various inequalities in the following forms or their reverses should hold with $\lambda_n=n^{\alpha}$ and $\Lambda_n$ defined as in \eqref{021} for different choices of $p$ and $\alpha$ with $U$ some constants depending on $p$ and $\lambda_n$'s whose values are suggested by the integral cases:
\begin{eqnarray}
\label{1.6}
   \sum^{\infty}_{n=1}\Big{(}\sum^{n}_{k=1} \frac {\lambda_ka_k}{\Lambda_n}
   \Big{)}^p & \leq & U \sum^{\infty}_{n=1}a_n^p, \\
  \sum^{\infty}_{n=1}\Big{(}\sum^{\infty}_{k=n} \frac {\lambda_ka_k}{\sum^{\infty}_{k=n} \lambda_k}
   \Big{)}^p & \leq & U \sum^{\infty}_{n=1}a_n^p, \nonumber \\
\label{1.8}
    \sum^{\infty}_{n=1}\frac {\Big{(}\sum^{\infty}_{k=n} \lambda_ka_k\Big{)}^p}{\Lambda^p_n}
    & \leq & U \sum^{\infty}_{n=1}a_n^p.
\end{eqnarray}

  For example, integral inequalities \eqref{1.2} suggest that inequalities \eqref{1.6} hold (for the special case $\lambda_n=n^{\alpha}$) for $p>1, (1+\alpha)p>1$ or $p<0, 1+\alpha>0$ with the reversed inequalities holding for $0<p<1, (1+\alpha)p>1$. The $p>1$ case has been extensively studied in \cite{G}, \cite{Be1}, \cite{G3} and will also be our main focus in the paper.
  We now take a look at the inequalities \eqref{1.8} for $0<p<1$. Integral inequalities \eqref{1.4} suggest that the reversed inequalities of \eqref{1.8} hold for the cases $0<p<1, 0<1+\alpha < 1/p$ (note that we get back \eqref{1.6'} for the special case $\alpha=0$ of $\lambda_n=n^{\alpha}$). We now treat these cases for general $\lambda_k$'s following the method in Section 3 of \cite{G3}. It is then easy to see that inequalities \eqref{1.8} hold for the cases $0<p<1, 0<1+\alpha < 1/p$ with $U=(p/(L-p))^p$ for $L>p$ if one can find a sequence ${\bf w}$ of positive terms with $w^{-1/(1-p)}_n\Lambda_n^{-p/(1-p)}$ decreasing to $0$, such
  that for any integer $n \geq 1$,
\begin{equation}
\label{3.1'}
 (w_1+\cdots+w_n)^{-1/(1-p)}\lambda^{-p/(1-p)}_n \leq
 \Big ( \frac {L-p}{p} \Big )^{p/(1-p)}
 \Big( \frac {w^{-1/(1-p)}_n}{\Lambda_n^{p/(1-p)}}- \frac {w^{-1/(1-p)}_{n+1}}{\Lambda_{n+1}^{p/(1-p)}}\Big).
\end{equation}
  We now define our sequence ${\bf w}$ inductively with $w_1=1$ and
\begin{equation*}
   \sum^{n}_{i=1}w_i=\Big ((1+\beta)\frac {\Lambda_n}{\lambda_n}-\beta \Big )w_n, \hspace{0.1in} \beta=\frac {2p-L}{L-p}.
\end{equation*}
   Note that this implies
\begin{equation*}
   w_{n+1}=\Big ((1+\beta)\frac {\Lambda_n}{\lambda_n}-\beta \Big )\frac {\lambda_{n+1}}{(1+\beta)\Lambda_n}w_n.
\end{equation*}
  We now set $x=\Lambda_n/\lambda_n$ and $y=\Lambda_{n+1}/\lambda_{n+1}$ and assume that $y \leq x+L$. Note that for our choice of ${\bf w}$, inequality \eqref{3.1'} can be recasted as
\begin{equation*}
   \Big ( 1+\frac {L/p-2}{x} \Big )^{1/(1-p)}-(y/x)^{1/(1-p)}(1-1/y)^{(1+p)/(1-p)} \geq \frac {L-p}{px}.
\end{equation*}
   It is easy to see the left-hand side expression above is a decreasing function of $y$ with $x$ fixed. Hence we may replace $y$ by $x+L$ and consider the following inequality for $x \geq 1$:
\begin{equation*}
   \Big ( 1+\frac {L/p-2}{x} \Big )^{1/(1-p)}-\Big(1+ \frac {L}{x} \Big )^{-p/(1-p)}\Big(1+ \frac {L-1}{x} \Big )^{(1+p)/(1-p)} \geq \frac {L-p}{px}.
\end{equation*}
   Suppose now we have $0<p \leq 1/3$, then by Taylor expansion, the left-hand side expression above is no less than
\begin{equation*}
   \Big ( 1+\frac {L/p-2}{(1-p)x} \Big )-\Big(1- \frac {pL}{(1-p)x} \Big )\Big(1+ \frac {(L-1)(1+p)}{(1-p)x} \Big ) \geq \frac {L-p}{px},
\end{equation*}
  provided that $L \geq 1$.

  Note that when $\lambda_k=k^{\alpha}$ with $-1<\alpha<0$, it follows from Lemma \ref{lem0} below (note that the right-hand side inequality of \eqref{6.1} still holds for $-1<r<0$) that we can take $L=1/(\alpha+1)$ here and it is then easy to check that
\begin{equation*}
   \frac {w^{-1/(1-p)}_n}{\Lambda_n^{p/(1-p)}} = O(n^{-(1-(2+\alpha)p+p^2(1+\alpha))/(p(1-p))}),
\end{equation*}
   so that $w_n^{-1/(1-p)}\Lambda_n^{-p/(1-p)}$ decreases to $0$ as $n$
   approaches infinity.  We now immediately deduce the following
\begin{theorem}
\label{thm3.2}
    Let $0<p \leq 1/3$ be fixed. The reversed inequalities \eqref{1.8} hold when $\lambda_n=n^{\alpha}$ for $-1< \alpha \leq 0$ with $U^{1/p} = (1+\alpha)p/(1-(1+\alpha)p)$.
\end{theorem}

   We point out here that even though one can often draw the analogues between the integral Hardy-type inequalities and the discrete ones, the two cases are
sometimes different. For example, in view of \eqref{1.6'} and the corresponding cases of \eqref{1.4}, one may suspect that inequalities \eqref{1.6'} hold for $0<p<1$. This is not the case, however, as one can take $a_1=1$ and $a_n=0$ in \eqref{1.6'} to see that the inequalities fail to hold when $(1-p)/p>1$ or $p>1/2$. In fact, it is also easy to show by setting $a_n=0$ for $n \geq 3$ that the case $p=1/2$ doesn't hold in \eqref{1.6'}.

   We end this section by pointing out other types of integral inequalities which are similar to Hardy-type inequalities. For the integral cases, we note the following inequalities that are similar to \eqref{0}:
\begin{eqnarray}
\label{2.11}
  \int^b_a s(x)|y(x)|^p|y'(x)|dx &\leq &  \int^b_ar(x)|y'(x)|^{p+1}dx, \hspace{0.1in} p>0, \\
\label{2.12}
  \int^b_a s(x)y^2(x)dx & \leq &  \int^b_ar(x)(y'(x))^{2}dx.
\end{eqnarray}
  where $r,s$ are non-negative measurable functions on $(a,b)$ and $y(x)$ absolutely continuous on $(a,b)$ subject to certain boundary conditions ($y(a)=0$ or $y(b)=0$ or both) with $r(x)|y'(x)|^{p+1}$ in \eqref{2.11}  and $r(x)(y'(x))^{2}$ in \eqref{2.12} integrable on $(a,b)$.

  The prototype of \eqref{2.11} is the case $a=0, p=1, r(x)=1, s(x)=4/b$ with $y(0)=y(b)=0$, which is referred to as the Opial's inequality in the literature and the prototype of \eqref{2.12} is the case $a=-\pi, b=\pi, r(x)=s(x)=1$ with $y(-\pi)=y(\pi), \int^{\pi}_{-\pi}y(x)dx=0$, which is referred to as the Wirtinger's inequality in the literature. We note here discrete analogues of Opial's inequality were studied by Wong \cite{W} and Lee \cite{Lee}, see also \cite{A&P} and the references therein for more details in this area. We also note that discrete analogues of Wirtinger's inequality were studied by Fan, Taussky and Todd \cite[Theorem 8]{FTT}, a subject we shall return to the discussion in Section \ref{sec 2'}.
\section{Carleman's approach versus Kaluza and Szeg\"o's approach}
\label{sec 2} \setcounter{equation}{0}
    We assume $p>1$ in this section and our goal in general is to find conditions on $\lambda_k$'s so that the following inequality holds
    for some constant $U$ and for any ${\bf a} \in l^p$:
\begin{equation*}
   \sum^{\infty}_{n=1}\Big{|}\sum^{n}_{k=1} \frac {\lambda_ka_k}{\Lambda_n}
   \Big{|}^p \leq U \sum^{\infty}_{n=1}|a_n|^p.
\end{equation*}
   It suffices to consider the cases with the infinite summations above replaced by any finite summations, say from $1$
   to $N \geq 1$ here. We may also assume $a_k \geq 0$ from now on and we shall define
\begin{equation*}
   A_n=\sum^n_{k=1}\frac {\lambda_ka_k}{\Lambda_n}.
\end{equation*}
   Carleman's approach is to determine the maximum value $\mu_N$ of $\sum^N_{n=1}A^p_n$ subject to the constraint
   $\sum^N_{n=1}a^p_n=1$ using Lagrange multipliers. We first show that we may further assume that $a_n > 0$ for all
   $1 \leq n \leq N$ when the maximum is reached. For otherwise, we may assume without loss of
   generality that $a_i=0, a_{i+1}>0$ for some $1 \leq i \leq N-1$ when the maximum is reached.
   We can now assume $a_n$'s are fixed for $n \neq i, i+1$.
   Then our assumption that $\sum^N_{n=1}a^p_n=1$ implies that the value of $a^p_i+a^p_{i+1}$ is constant and hence defines
   $a_{i+1}$ explicitly as a function of $a_i$.
   We now regard $\sum^N_{n=1}A^p_n$ as a function of $a_i$ and it is then easy to check that it is an increasing function of $a_i$ near $a_i=0$,
   by which it means that on increasing the value of $a_i$ from $0$ to a small
   positive number while decreasing the value of $a_{i+1}$ and keeping other variables as well as the sum $\sum^N_{n=1}a^p_n=1$ fixed,
   we will increase the value of $\sum^N_{n=1}A^p_n$, a contradiction.

  We now define
\begin{equation*}
  F({\bf a}; \mu)=\sum^N_{n=1}A^p_n-\mu (\sum^N_{n=1}a^p_n-1),
\end{equation*}
  where ${\bf a}=(a_n)_{1 \leq n \leq N}$. By the Lagrange method and our discussions above, we have to solve $\nabla F=0$, or the following system of equations:
\begin{equation}
\label{2.1}
  \mu a^{p}_k=\sum^N_{n=k}\frac {\lambda_kA^{p-1}_n}{\Lambda_n}a_k, \hspace{0.1in} 1 \leq k \leq N; \hspace{0.1in} \sum^N_{n=1}a^p_n=1.
\end{equation}
   We note that on summing over $1 \leq k \leq N$ of the first $N$ equations above, we get
\begin{equation*}
  \sum^N_{n=1}A^p_n=\mu.
\end{equation*}
   Hence we have $\mu=\mu_N$ in this case, which allows us to recast the equations \eqref{2.1} as:
\begin{equation}
\label{2.10}
  \mu_N \frac {a^{p-1}_k}{\lambda_k}=\sum^N_{n=k}\frac {A^{p-1}_n}{\Lambda_n}, \hspace{0.1in} 1 \leq k \leq N; \hspace{0.1in} \sum^N_{n=1}a^p_n=1.
\end{equation}
  On subtracting consecutive equations, we can rewrite the above system of equations as:
\begin{equation*}
  \mu_N (\frac {a^{p-1}_k}{\lambda_k}-\frac {a^{p-1}_{k+1}}{\lambda_{k+1}})=\frac {A^{p-1}_k}{\Lambda_k}, \hspace{0.1in} 1 \leq k \leq N-1; \hspace{0.1in}  \mu_N \frac {a^{p-1}_N}{\lambda_N}=\frac {A^{p-1}_N}{\Lambda_N}; \hspace{0.1in} \sum^N_{n=1}a^p_n=1.
\end{equation*}

   Now we define for $1 \leq k \leq N-1$,
\begin{equation*}
  \omega_k = \frac {\Lambda_k}{\lambda_k}-\frac {\Lambda_k a^{p-1}_{k+1}}{\lambda_{k+1}a^{p-1}_k},
\end{equation*}
  so that we can further rewrite our system of equations as:
\begin{equation*}
  \mu_N a^{p-1}_k \omega_k=A^{p-1}_k, \hspace{0.1in} 1 \leq k \leq N-1; \hspace{0.1in}  \mu_N \frac {a^{p-1}_N}{\lambda_N}=\frac {A^{p-1}_N}{\Lambda_N}; \hspace{0.1in} \sum^N_{n=1}a^p_n=1.
\end{equation*}
   It is easy to check that for $1 \leq k \leq N-2$,
\begin{equation*}
   \omega^{\frac 1{p-1}}_{k+1}=\frac {\Lambda_{k}}{\Lambda_{k+1}}\Big (\frac {\omega_{k}}{\frac {\lambda_{k+1}}{\Lambda_k}(\Lambda_k/\lambda_k-\omega_{k})} \Big )^{\frac 1{p-1}}+\frac {\lambda_{k+1}}{\Lambda_{k+1}} \Big ( \frac 1{ \mu_N} \Big )^{\frac 1{p-1}}.
\end{equation*}
   We now define a sequence of real functions $\Omega_k(\mu)$ inductively by setting $\Omega_1(\mu)=1/\mu$ and
\begin{equation*}
   \Omega^{\frac 1{p-1}}_{k+1}(\mu)=\frac {\Lambda_{k}}{\Lambda_{k+1}}\Big (\frac {\Omega_{k}(\mu)}{\frac {\lambda_{k+1}}{\Lambda_k}(\Lambda_k/\lambda_k-\Omega_{k}(\mu))} \Big )^{\frac 1{p-1}}+\frac {\lambda_{k+1}}{\Lambda_{k+1}} \Big ( \frac 1{ \mu} \Big )^{\frac 1{p-1}}.
\end{equation*}
   We note that $\Omega_k(\mu_N)= \omega_k$ for $1 \leq k \leq N-1$ and
\begin{eqnarray*}
   \Omega^{\frac 1{p-1}}_N(\mu_N) &=& \frac {\Lambda_{N-1}}{\Lambda_{N}}\Big (\frac {\omega_{N-1}}{\frac {\lambda_{N}}{\Lambda_{N-1}}(\Lambda_{N-1}/\lambda_{N-1}-\omega_{N-1})} \Big )^{\frac 1{p-1}}+\frac {\lambda_{N}}{\Lambda_{N}} \Big ( \frac 1{ \mu_N} \Big )^{\frac 1{p-1}} \\
 &=&  \frac {\Lambda_{N-1}}{\Lambda_{N}}\Big (\frac {A^{p-1}_{N-1}}{\mu_Na^{p-1}_N} \Big )^{\frac 1{p-1}}+\frac {\lambda_{N}}{\Lambda_{N}} \Big ( \frac 1{ \mu_N} \Big )^{\frac 1{p-1}}=\Big ( \frac 1{ \mu_N} \Big )^{\frac 1{p-1}}\frac {A_N}{a_N} \\
 &=& \Big ( \frac {\Lambda_N}{ \lambda_N} \Big )^{\frac 1{p-1}} .
\end{eqnarray*}
   We now define another sequence of real functions $\eta_k(\mu)$ by setting
\begin{equation*}
  \eta_k(\mu)=\Big (\frac {\Lambda_k}{\lambda_k}\Big)^{p-1} \Omega_k(\mu)
\end{equation*}
 so that it satisfies the following relations:
\begin{equation}
\label{2.2'}
   \eta^{\frac 1{p-1}}_{k+1}(\mu)=\frac {\Lambda_{k}}{\lambda_{k+1}}\Big (\frac {\Lambda_k\eta_{k}(\mu)/\lambda_{k+1}}{(\Lambda_k/\lambda_k)^p-\eta_{k}(\mu)} \Big )^{\frac 1{p-1}}+ \Big ( \frac 1{ \mu} \Big )^{\frac 1{p-1}}.
\end{equation}
    Note that we have seen above that $\eta_N(\mu_N)=(\Lambda_{N}/\lambda_{N})^p$ and Carleman's idea is to show that
    the above relations \eqref{2.2'} lead to a contradiction if $\mu$ is large and this forces  $\mu_N$ to be small.
    For example, one can show by induction that if \eqref{024} is satisfied and $\mu > (1-L/p)^{-p}$, then for  $k \geq 1$,
\begin{equation*}
   \eta^{\frac 1{p-1}}_k(\mu) < (b+c)\Big ( \frac {\Lambda_k}{\lambda_k}\Big )-c, \hspace{0.1in} b=(1-L/p)^{p/(p-1)}, \hspace{0.1in}  c=\frac {L}{p}(1-L/p)^{1/(p-1)}.
\end{equation*}
   It follows that if the above assertion is established, then for $1 \leq n \leq N$,
\begin{equation*}
  0< \eta^{\frac 1{p-1}}_n(\mu) <(b+c)\Big (\frac {\Lambda_n}{\lambda_n} \Big )-c < (b+c)\frac {\Lambda_n}{\lambda_n} < \frac {\Lambda_n}{\lambda_n} \leq \Big ( \frac
  {\Lambda_n}{\lambda_n} \Big )^{\frac p{p-1}}.
\end{equation*}
   As we have seen above that $\eta_N(\mu_N)=(\Lambda_{N}/\lambda_{N})^p$, this forces $\mu_N \leq (1-L/p)^{-p}$ and the assertion
   of Theorem \ref{thm03} will follow.  We further note here that one can compare the above approach
   to the case considered in Section \ref{sec 2'}, where the $l^2$ norms of weighted mean matrices is treated first using linear algebra techniques.
   Then the general $l^p$ cases will be treated using the duality principles. It is easy to see that the method used there can be regarded
   essentially as the special case $\mu=(p/(p-L))^p$ in the outlined proof above of Theorem \ref{thm03}. We shall leave the details to Section \ref{sec 2'}.

   We now give a short account on Kaluza and Szeg\"o's approach \cite{K} on proving Hardy's inequality \eqref{eq:1}. In fact,
   we shall explain this for the general case involving weighted mean matrices. Using the notations
in Section \ref{sec 1} and once again restricting our attention to
any finite summations, say from $1$ to $N \geq 1$ here, we are
  looking for a positive constant $U$ such that
\begin{equation}
\label{2.0}
   \sum^{N}_{n=1}\Big{|} \frac {1}{\Lambda_n}\sum^{n}_{k=1}\lambda_ka_k
   \Big{|}^p \leq U \sum^{N}_{n=1}|a_n|^p
\end{equation}
   holds for all complex sequences ${\bf a} \in l^p$. To motivate the approach, we may assume $a_n \geq 0$ \
   and we are using Carleman's approach to find the maximum value $\mu_N$ of $\sum^N_{n=1}A^p_n$ subject to the constraint $\sum^N_{n=1}a^p_n=1$.
   Suppose this is done and we find that the maximum value is reached at a sequence ${\bf w}=(w_n)_{1 \leq n \leq N}$. Hence \eqref{2.1} is satisfied
   with $a_n$'s there replaced by $w_n$'s and $\mu=\mu_N$. This motivates us to consider,
   for an arbitrary sequence ${\bf a}=(a_n)_{1 \leq n \leq N}$, the following expression
\begin{equation*}  \mu_N\sum^N_{n=1}|a_n|^p=\sum_{k=1}^{N}w_k^{-(p-1)}\Big(\sum_{n=k}^{N}\frac
  {\lambda_k}{\Lambda_n}\Big(\sum_{j=1}^n\frac
  {\lambda_j}{\Lambda_n}w_j\Big)^{p-1}\Big)|a_k|^p.
\end{equation*}
  Thus inequality \eqref{2.0} will follow from this with $U=\mu_N$ if one can show the right-hand side expression above is no less than the left-hand side expression of \eqref{2.0}.

   Kaluza and Szeg\"o's idea is to reverse the process discussed above by finding an
  auxiliary sequence ${\bf w}=(w_n)_{1 \leq n \leq N}$ of positive terms
  such that by H\"older's inequality,
\begin{eqnarray*}
\label{eq:12}
  \Big( \sum_{k=1}^n\lambda_k|a_k| \Big)^p &= & \Big( \sum_{k=1}^n\lambda_k|a_k|w_k^{-\frac {1}{q}} \cdot
  w_k^{\frac {1}{q}} \Big )^p \\
  & \leq & \Big( \sum_{k=1}^n\lambda^p_k|a_k|^pw_k^{-(p-1)} \Big ) \Big ( \sum_{j=1}^n w_j \Big)^{p-1}
\end{eqnarray*}
  so that
\begin{eqnarray*}
\label{eq:13}
 \sum^{N}_{n=1}\Big{|} \frac {1}{\Lambda_n}\sum^{n}_{k=1}\lambda_ka_k
   \Big{|}^p  &\leq &
 \sum^N_{n=1}\frac {1}{\Lambda^p_n}\Big(\sum_{k=1}^{n}\lambda^p_k|a_k|^pw_k^{-(p-1)}\Big)\Big(\sum_{j=1}^n
 w_j\Big)^{p-1} \\
  &=& \sum_{k=1}^{N}w_k^{-(p-1)}\lambda^p_k\Big(\sum_{n=k}^{N}\frac
  {1}{\Lambda^p_n}\Big(\sum_{j=1}^nw_j\Big)^{p-1}\Big)|a_k|^p.
\end{eqnarray*}
  An ideal approach would now for one to be able to find for each $p>1$ a positive constant $U$,
  a sequence ${\bf w}$ of positive terms, such
  that for any integer $1 \leq n \leq N$,
\begin{equation}
\label{3.70}
  \lambda^p_k\Big(\sum_{n=k}^{N}\frac
  {1}{\Lambda^p_n}\Big(\sum_{j=1}^nw_j\Big)^{p-1}\Big)=Uw_k^{p-1}.
\end{equation}
   One checks easily that the above relations give back the first $N$ relations in \eqref{2.10} upon a change of variables $w_k \rightarrow \lambda_ka_k$ and this will certainly make inequality \eqref{2.0} hold. Moreover, in this case, it is easy to see that when $\lambda_k \neq 0$, the optimal $a_k$'s appearing in \eqref{2.0} will take the values $|a_k|=c \cdot w^{2/q}_k/\lambda_k$ for some positive constant $c$.

  Since it is difficult in general to find a sequence ${\bf w}$ to satisfy the conditions \eqref{3.70}, one may suppose now that one can find for each $p>1$ a positive constant $U$,
  a sequence ${\bf w}$ of positive terms, such
  that for any integer $1 \leq n \leq N$,
\begin{equation}
\label{eq:7}
 \Big(\sum_{i=1}^nw_i \Big )^{p-1} \leq U\Lambda^p_n \Big ( \frac {w_n^{p-1}}{\lambda^p_n}-\frac {w_{n+1}^{p-1}}{\lambda^p_{n+1}} \Big ),
\end{equation}
  where we define $w_{N+1}=0$. Then it is easy to see that inequality \eqref{2.0} follows from this.
  When $\lambda_n=1$ for all $n$, Kaluza and Szeg\"o's choice for ${\bf w}$ is given inductively by setting $w_1=1$ and
\begin{equation*}
   \sum_{i=1}^nw_i=\frac {n-1/p}{1-1/p}w_n.
\end{equation*}
 and one can show that \eqref{eq:7}
  holds in this case with $U=q^p$ and Hardy's inequality \eqref{eq:1}
  follows from this.

\section{Proof of Theorem \ref{thm03}}
\label{sec 3} \setcounter{equation}{0}
       We now apply Kaluza and Szeg\"o's method to give a proof of Theorem \ref{thm03}.  We note first by a change of variables $w_k \rightarrow \lambda_ka_k$, we can recast \eqref{eq:7} as
\begin{equation}
\label{4.3}
  U \Big (\frac {a^{p-1}_k}{\lambda_k}-\frac {a^{p-1}_{k+1}}{\lambda_{k+1}} \Big ) \geq \frac {A^{p-1}_k}{\Lambda_k}, \hspace{0.1in} 1 \leq k \leq N-1; \hspace{0.1in}  U \frac {a^{p-1}_N}{\lambda_N} \geq \frac {A^{p-1}_N}{\Lambda_N},
\end{equation}
  where $A_n$'s are defined as in the previous section. It now suffices to find a sequence ${\bf a}=(a_n)_{1 \leq n \leq N}$ of positive terms so that inequalities \eqref{4.3} are satisfied with $U=(p/(p-L))^p$. We now define our sequence inductively by setting $a_1=1$ and for $n \geq 1$,
  \begin{equation*}
     A_n=\frac 1{\Lambda_n}\sum^n_{i=1}\lambda_ia_i=(1+\beta-\frac {\beta \lambda_n}{\Lambda_n})a_n,
  \end{equation*}
    where $\beta=L/(p-L)$. Equivalently, this is amount to taking $\sum^n_{i=1}w_i=((1+\beta)\Lambda_n/\lambda_n-\beta)w_n$ for those $w_i$'s satisfying \eqref{eq:7}. One checks easily that the above relations lead to
  the following relation between $a_n$ and $a_{n+1}$:
  \begin{equation*}
     a_{n+1}=\frac 1{1+\beta}(1+\beta-\frac {\beta \lambda_n}{\Lambda_n})a_n.
  \end{equation*}
      It is then easy to see that inequalities \eqref{4.3} follow from the following inequality for $n \geq 1$:
  \begin{equation*}
     U\Big ( \frac {\Lambda_n}{\lambda_n} (1+\beta-\frac {\beta \lambda_n}{\Lambda_n})^{1-p}-(1+\beta)^{1-p}(\frac {\Lambda_{n+1}}{\lambda_{n+1}}-1)\Big ) \geq 1.
  \end{equation*}
      We now set $x=\Lambda_n/\lambda_n$, $y= \Lambda_{n+1}/\lambda_{n+1}$ to rewrite the above inequality as:
 \begin{equation*}
     U(1+\beta)^{1-p}\Big ( x\Big (1-\frac {\beta }{(1+\beta)x} \Big )^{1-p}-(y-1)\Big ) \geq 1.
 \end{equation*}
     The above inequality now follows from \eqref{024} and this completes the proof of Theorem \ref{thm03}.

\section{Duality and Another Proof of Theorem \ref{thm03}}
\label{sec 2'} \setcounter{equation}{0}
   We first recall the following duality principle concerning the norms of linear operators:
\begin{theorem}\cite[Lemma 2]{M}
\label{thm5} Let $p > 1$ be fixed and let $C=(c_{n,k})$ be a fixed $N \times K$ matrix. Then the following three assertions concerning the constant $U$ for any ${\bf x} \in l^p, {\bf y} \in l^q$ are equivalent:
\begin{eqnarray*}
 && \Big ( \sum_{n} \Big | \sum_{k} c_{n,k}x_k \Big |^p \Big )^{1/p} \leq  U ||{\bf x}||_p, \\
  && \Big ( \sum_{k} \Big | \sum_{n} c_{n,k}y_n \Big |^q \Big )^{1/q} \leq  U ||{\bf y}||_q, \\
 && \Big | \sum_{n,k} c_{n,k}x_ky_n \Big | \leq  U ||{\bf x}||_p||{\bf y}||_q.
\end{eqnarray*}
\end{theorem}

  We now describe Wang and Yuan's method in \cite{W&Y} for establishing Corollary \ref{cor03} for $p=2$. We may assume $a_n$ being real without loss of generality and
   it suffices to prove the corollary for any finite summations from $n=1$ to $N$ with $N \geq 1$.
   We also note that it follows from our assumption on $L$ that $\lambda_n>0$. Now consider
\begin{equation*}
   \sum^{N}_{n=1}\Big ( \sum^{n}_{i=1}\frac {\lambda_i}{\Lambda_n}a_i\Big )^2
   =\sum^{N}_{n=1}\Big ( \sum^{n}_{i,j=1}\frac {\lambda_i\lambda_j}{\Lambda^2_n}a_ia_j\Big )
   =\sum^{N}_{i,j=1}\alpha_{i,j}a_ia_j, \hspace{0.1in}
   \alpha_{i,j}= \sum^{N}_{k \geq \max{(i,j)}}\frac {\lambda_i\lambda_j}{\Lambda^2_k}.
\end{equation*}
  We view the above as a quadratic form and define the associated
  matrix $A$ to be
\begin{equation}
\label{4.1}
  A=\Big ( \alpha_{i, j} \Big )_{1 \leq i, j \leq N}.
\end{equation}
   We note that the matrix $A$ here is certainly positive definite, being equal to $B^{t}B$
   with $B$ a lower-triangular matrix,
\begin{equation*}
   B=\Big ( b_{i, j} \Big )_{1 \leq i, j \leq N},  \hspace{0.1in} b_{i,j}=\lambda_j/\Lambda_i,  ~~ 1 \leq j \leq
    i;   \hspace{0.1in} b_{i,j}=0, \hspace{0.1in} j > i.
\end{equation*}
   It is easy to check that the entries of $B^{-1}$ are given by
\begin{equation*}
  \big (  B^{-1} \big )_{i,i}=\frac
  {\Lambda_i}{\lambda_i}, \hspace{0.1in}
  \big (  B^{-1} \big )_{i+1,i}=-\frac
  {\Lambda_i}{\lambda_{i+1}}, \hspace{0.1in} \big (  B^{-1} \big
  )_{i,j}=0 \hspace{0.1in} \text{otherwise}.
\end{equation*}

    In order to
   establish our assertion, it suffices to show that the maximum eigenvalue of $A$ is less than
   $4/(2-L)^2$ or the minimum eigenvalue of its inverse $A^{-1}$ is greater than
   $(2-L)^2/4$ and this is equivalent to proving that the matrix $A^{-1} - \lambda I_N$ is positive
   definite, where $\lambda=(2-L)^2/4$ and $I_N$ is the $N \times N$ identity matrix. Using the expression $A^{-1}=B^{-1}(B^{-1})^{t}$, we
   see that this is equivalent to showing that
   for any integer $N \geq 1$ and any
   real sequence ${\bf a}=(a_n)_{1 \leq n \leq N}$,
\begin{equation}
\label{3.0}
   \sum^{N-1}_{n=1}\Big ( \frac
  {\Lambda_n}{\lambda_n}a_n-\frac
  {\Lambda_n}{\lambda_{n+1}}a_{n+1}\Big )^2 +  \frac
  {\Lambda^2_N}{\lambda^2_N}a^2_N\geq \frac {(2-L)^2}{4}\sum^N_{n=1}a_n^2.
\end{equation}

  We point out here that this linear algebraic approach can be viewed as an approach via the duality principle of linear operators. In our case, we may regard the matrix $A$ given in \eqref{4.1} as a linear operator on the $l^2$ space, which is self-dual. Hence its adjoint $A^t$ also acts on $l^2$ with the same norm as $A$ by Theorem \ref{thm5}. If we reformulate this in terms of inequalities, what we need to prove are the following Copson (see \cite[Theorems 331, 344]{HLP}) type inequalities:
\begin{equation*}
  \sum^{N}_{n=1}\Big ( \sum^{N}_{k=n}\frac {\lambda_n}{\Lambda_k}a_k\Big )^2 \leq \frac {4}{(2-L)^2}\sum^N_{n=1}a^2_n.
\end{equation*}
   Now by a change of variables $\sum^{N}_{k=n}\frac {\lambda_n}{\Lambda_k}a_k=b_n$, one sees easily that the above inequalities give back \eqref{3.0} (identifying $b_n$ with $a_n$ there). In the special case of $\lambda_n=1$, Wang and Yuan \cite[Theorem 2]{W&Y} obtained a refinement of Hardy's inequality for $p=2$ via their approach and we note here that if one uses the duality approach, then one can obtain the same result, saving the effort of inverting a matrix.

  For any integer $n \geq 1$ and fixed constants $\alpha$, $\beta$, $a_{n+1}, \mu_n$ (here $\alpha$, $\beta$ may depend on $n$), we consider the following function:
\begin{equation*}
  f(a_n) :=(\alpha a_n-\beta a_{n+1})^2-\mu_na^2_n.
\end{equation*}
   When $\mu_n > \alpha^2$, it is easy to see that
\begin{equation}
\label{2.2}
  f(a_n) \leq f(\frac {\alpha \beta a_{n+1}}{\alpha^2-\mu_n})=\frac {\beta^2\mu_{n}\alpha^2_{n+1}}{\mu_n-\alpha^2},
\end{equation}
  with the above inequality reversed when $\mu_n < \alpha^2$.

  On taking $\alpha=\Lambda_n/\lambda_n, \beta=\Lambda_n/\lambda_{n+1}$ here, we obtain that for any $0 < \mu_n < \Lambda^2_n/\lambda^2_n$,
\begin{equation*}
   \Big ( \frac
  {\Lambda_n}{\lambda_n}a_n-\frac
  {\Lambda_n}{\lambda_{n+1}}a_{n+1}\Big )^2- \mu_n a^2_n \geq -\frac {\Lambda^2_n/\lambda^2_{n+1}}
  {\Lambda^2_n/\lambda^2_n-\mu_n}\mu_na^2_{n+1}.
\end{equation*}
   Summing the above inequality for $n=1, \ldots, N-1$ yields:
\begin{eqnarray*}
 &&  \sum^{N-1}_{n=1}\Big ( \frac
  {\Lambda_n}{\lambda_n}a_n-\frac
  {\Lambda_n}{\lambda_{n+1}}a_{n+1}\Big )^2 +  \frac
  {\Lambda^2_N}{\lambda^2_N}a^2_N \\
  &\geq &  \mu_1a^2_1+\sum^{N-2}_{n=1}\Big(\mu_{n+1}- \frac {\Lambda^2_n/\lambda^2_{n+1}}
  {\Lambda^2_n/\lambda^2_n-\mu_n}\mu_n\Big )a_{n+1}^2+ \Big ( \frac
  {\Lambda^2_N}{\lambda^2_N}  - \frac {\Lambda^2_{N-1}/\lambda^2_N}
  {\Lambda^2_{N-1}/\lambda^2_{N-1}-\mu_{N-1}}\mu_{N-1} \Big ) a^2_N.
  \nonumber
\end{eqnarray*}
   Now it suffices to show that one can always find a sequence of numbers $(\mu_n)_{n \geq 1}$ with $\mu_1=(2-L)^2/4, 0 < \mu_n < \Lambda^2_n/\lambda^2_n$ for any integer $n \geq 1$,
   such that the following relations hold for $1 \leq n \leq N-1$:
\begin{equation*}
  \mu_{n+1}- \frac {\Lambda^2_n/\lambda^2_{n+1}}
  {\Lambda^2_n/\lambda^2_n-\mu_n}\mu_n = \frac {(2-L)^2}{4}.
\end{equation*}
   Note that this is just a special case of \eqref{2.2'} of $p=2$, $\mu=4/(2-L)^2$ if we identify the sequence $(\mu_n)_{n \geq 1}$ with
   the sequence $(\eta_n)_{n \geq 1}$ there.


    We now consider the approach via the duality principle for the general case, which allows us to give another proof of Theorem \ref{thm03}.
    Let $p>1$ be fixed and let $N \geq 1$ be an integer, we now seek for conditions on the $\lambda_k$'s
    such that inequality \eqref{2.0} holds for any ${\bf a} \in
    l^p$ with $U=(p/(p-L))^p$ for some $0<L<p$. By the duality principle Theorem \ref{thm5},
    this is equivalent to establishing the following inequality
\begin{equation}
\label{2.100}
   \sum^{N}_{n=1}\Big{|}\sum^{N}_{k=n} \frac {\lambda_na_k}{\Lambda_k}
   \Big{|}^q \leq \Big (\frac {p}{p-L} \Big )^q \sum^{N}_{n=1}|a_n|^q.
\end{equation}
  Without loss of generality, we may assume that all the $a_k$'s
  are non-negative and define
\begin{equation*}
    \sum^{N}_{k=n} \frac {\lambda_na_k}{\Lambda_k}=b_n.
\end{equation*}
   This allows us to recast inequality \eqref{2.100} as
\begin{equation}
\label{2.20}
   \sum^{N}_{n=1}b^q_n \leq \Big (\frac {p}{p-L} \Big )^q \sum^{N}_{n=1}
   \Big(\frac {\Lambda_n}{\lambda_n}b_n-\frac {\Lambda_n}{\lambda_{n+1}}b_{n+1}\Big
   )^q,
\end{equation}
  where we set $b_{N+1}=0$ and assume that $b_k \geq 0,
  b_k/\lambda_k \geq b_{k+1}/\lambda_{k+1}, 1 \leq k \leq N-1$.

  Fixing a set of positive $\lambda_k$'s as well as $b_{n+1}$ and choose a non-negative number $\mu_n$,
  whose value is to be determined later, we consider now for $x \geq
  \lambda_nb_{n+1}/\lambda_{n+1}$,
\begin{equation*}
   f(x)=\Big(\frac {\Lambda_n}{\lambda_n}x-\frac
{\Lambda_n}{\lambda_{n+1}}b_{n+1}\Big
   )^q-\mu_nx^q.
\end{equation*}
   It is readily checked that $f'(x)=0$ implies that
\begin{equation}
\label{2.30}
    \frac {\Lambda_n}{\lambda_n}x-\frac
{\Lambda_n}{\lambda_{n+1}}b_{n+1}=\Big (\frac
{\mu_n\lambda_n}{\Lambda_n} \Big )^{1/(q-1)}x.
\end{equation}
  Solving this for $x$, we obtain
\begin{equation*}
    x=\frac {\Lambda_nb_{n+1}/\lambda_{n+1}}{\Lambda_n/\lambda_n-\Big (\mu_n\lambda_n/\Lambda_n \Big )^{1/(q-1)}}.
\end{equation*}
  It follows that $x \geq 0$ implies that
  $\mu_n<(\Lambda_n/\lambda_n)^q$, a condition we shall enforce
  from now on. With this condition on $\mu_n$, it is then easy to
  check that for $x$ taking the above value, we have $f''(x) \geq
  0$ and it follows that (note that $x/\lambda_n \geq
  b_{n+1}/\lambda_{n+1}$ thanks to \eqref{2.30})
\begin{eqnarray*}
   f(b_n) & \geq & \frac {\Big (\frac
{\mu_n\lambda_n}{\Lambda_n} \Big )^{q/(q-1)}-\mu_n}{\Big ( \frac
{\Lambda_n}{\lambda_{n}}-\Big (\frac {\mu_n\lambda_n}{\Lambda_n}
\Big )^{1/(q-1)}\Big )^q}\Big ( \frac
{\Lambda_n}{\lambda_{n+1}}\Big )^qb^q_{n+1} \\
 &=& -\frac {1}{\Big
(\Big ( \Lambda_n/\lambda_{n} \Big )^{q/(q-1)}\mu_n^{-1/(q-1)}-1
\Big )^{q-1}}\Big ( \frac {\Lambda_n}{\lambda_{n+1}}\Big
)^qb^q_{n+1}.
\end{eqnarray*}
   We now recast the above inequality as
\begin{equation*}
   \Big(\frac {\Lambda_n}{\lambda_n}b_n-\frac
{\Lambda_n}{\lambda_{n+1}}b_{n+1}\Big
   )^q \geq \mu_nb^q_n- \frac {1}{\Big (\Big ( \Lambda_n/\lambda_{n} \Big
)^{q/(q-1)}\mu_n^{-1/(q-1)}-1 \Big )^{q-1}}\Big ( \frac
{\Lambda_n}{\lambda_{n+1}}\Big )^qb^q_{n+1},
\end{equation*}
  and on adding up both sides above for $1 \leq n \leq N$, we obtain
\begin{eqnarray*}
  && \sum^{N}_{n=1}
   \Big(\frac {\Lambda_n}{\lambda_n}b_n-\frac {\Lambda_n}{\lambda_{n+1}}b_{n+1}\Big
   )^q \\
 &\geq & \mu_1b^q_1+\sum^{N-1}_{n=1}\Big(\mu_{n+1}-\frac {1}{\Big
(\Big ( \Lambda_n/\lambda_{n} \Big )^{q/(q-1)}\mu_n^{-1/(q-1)}-1
\Big )^{q-1}}\Big ( \frac {\Lambda_n}{\lambda_{n+1}}\Big )^q \Big
)b^q_{n+1}.
\end{eqnarray*}

   We now set
\begin{equation*}
 \mu_1=\mu_{n+1}-\frac {1}{\Big
(\Big ( \Lambda_n/\lambda_{n} \Big )^{q/(q-1)}\mu_n^{-1/(q-1)}-1
\Big )^{q-1}}\Big ( \frac {\Lambda_n}{\lambda_{n+1}}\Big )^q =\Big
( \frac {p-L}{p} \Big )^q:=b.
\end{equation*}
  This will lead to \eqref{2.20} provided the condition
  $\mu_n<(\Lambda_n/\lambda_n)^q$ is satisfied. Note that this corresponds to the case of \eqref{2.2'} for $\mu=(p/(p-L))^p$ if we identify
  the sequence $(\mu_n)_{n \geq 1}$ with
   the sequence $(\eta^{1/(p-1)}_n)_{n \geq 1}$ there. We now proceed inductively to see what conditions
  will be imposed on the $\lambda_n$'s so that we can have $\mu_n \leq (b+c)\Lambda_n/\lambda_n-c$ with $c$ a constant
  to be specified later so that $\mu_n<(\Lambda_n/\lambda_n)^q$ is
  satisfied. First note that when
  $n=1$,  our choice of $\mu_1$ posts no restrictions on the
  $\lambda_n$'s. Suppose now $\mu_n \leq
  (b+c)\Lambda_n/\lambda_n-c$ is satisfied and we then want to
  show
\begin{eqnarray*}
   \mu_{n+1} &=& \frac {1}{\Big
(\Big ( \Lambda_n/\lambda_{n} \Big )^{q/(q-1)}\mu_n^{-1/(q-1)}-1
\Big )^{q-1}}\Big ( \frac {\Lambda_n}{\lambda_{n+1}}\Big )^q +\Big
( \frac {p-L}{p} \Big )^q \\
 &\leq &
(b+c)\Lambda_{n+1}/\lambda_{n+1}-c=(b+c)\Lambda_{n}/\lambda_{n+1}+b.
\end{eqnarray*}
   One is then led to show that
\begin{equation*}
 (b+c)^{-1/(q-1)}\frac {\Lambda_n}{\lambda_{n+1}} \leq \Big ( \Lambda_n/\lambda_{n} \Big
 )^{q/(q-1)}\mu_n^{-1/(q-1)}-1.
\end{equation*}
  Using our assumption on $\mu_n$, it then remains to show that
\begin{equation*}
 (b+c)^{-1/(q-1)}\frac {\Lambda_n}{\lambda_{n+1}} \leq \Big ( \Lambda_n/\lambda_{n} \Big
 )^{q/(q-1)}\Big((b+c)\Lambda_n/\lambda_n-c \Big )^{-1/(q-1)}-1.
\end{equation*}
  We recast the above inequality in terms of $p$ to get
\begin{equation}
\label{2.4}
   \frac {\Lambda_n}{\lambda_{n+1}} \leq \frac {\Lambda_n}{\lambda_n}
 \Big(1-\frac {c}{b+c}\frac {\lambda_n}{\Lambda_n} \Big )^{1-p}-(b+c)^{p-1}.
\end{equation}
   This is the condition to be satisfied by the $\lambda_n$'s. Now to determine the value of $c$, we apply the Taylor expansion of the
first term of the
   right-hand side expression above to see that we need to choose
   $c$ so that the expression
\begin{equation*}
   \frac {(p-1)c}{b+c}-(b+c)^{p-1}
\end{equation*}
   is maximized. It is then easy to check that we need to take
\begin{equation*}
   c=\Big (\frac {p-L}{p} \Big )^{1/(p-1)} \cdot \frac {L}{p}.
\end{equation*}
   It is readily seen that $\mu_n \leq (b+c)\Lambda_n/\lambda_n-c<(\Lambda_n/\lambda_n)^q$ is satisfied and \eqref{2.4} becomes \eqref{024}.
   It follows that
   inequality \eqref{2.100} holds as long as \eqref{024} is satisfied for $1 \leq n \leq N$. As $N$ is arbitrary, this now concludes our proof of
   Theorem \ref{thm03}.

   We point out here inequality \eqref{3.0} can be regarded as an analogue to the following discrete inequality of Wirtinger's type studied by Fan, Taussky and Todd \cite[Theorem 8]{FTT}:
\begin{equation}
\label{3.11}
  a^2_1+\sum^{N-1}_{n=1}(a_n-a_{n+1})^2 + a^2_N \geq 2 \Big (1-\cos \frac {\pi}{N+1} \Big )\sum^{N}_{n=1} a^2_n.
\end{equation}
    The converse of the above inequality was found by I. \v Z. Milovanovi\'c and
    G. V. Milovanovi\'c \cite{MM}:
\begin{equation}
\label{3.22}
  a^2_1+\sum^{N-1}_{n=1}(a_n-a_{n+1})^2 + a^2_N \leq 2 \Big (1+\cos \frac {\pi}{N+1} \Big )\sum^{N}_{n=1} a^2_n.
\end{equation}
   Simple proofs of inequalities \eqref{3.11} and \eqref{3.22} were given by Redheffer \cite{R2} and Alzer \cite{A1}, respectively. Our proof of Theorem \ref{thm03} given above is motivated by the methods used in  \cite{R2} and \cite{A1}.

  To end this section, we note the paper \cite{Lo} contains several generalizations of inequalities \eqref{3.11} and \eqref{3.22}, one of them can be stated as:
\begin{theorem}
\label{thm3.1}
  For any real sequence ${\bf a}=(a_n)_{1 \leq n \leq N}$, and two positive real numbers $a, b$,
\begin{equation}
\label{3.7}
  (a^2+b^2-2ab\cos \frac {\pi}{N+1} \Big )\sum^{N}_{n=1} a^2_n \leq b^2a^2_1+\sum^{N-1}_{n=1}(aa_n-ba_{n+1})^2 + a^2a^2_N \leq  \Big (a^2+b^2+2ab\cos \frac {\pi}{N+1} \Big )\sum^{N}_{n=1} a^2_n.
\end{equation}
\end{theorem}

  The proof given in \cite{Lo} to the above theorem is to regard
\begin{equation*}
  b^2a^2_1+\sum^{N-1}_{n=1}(aa_n-ba_{n+1})^2 + a^2a^2_N
\end{equation*}
   as a quadratic form with the associated
  matrix $A$ being symmetric tridiagonal whose entries are
   given by
\begin{equation*}
  \big (  A \big )_{i,i}=a^2+b^2, \hspace{0.1in}
  \big (  A \big )_{i,i+1}=\big ( A \big )_{i+1,i}=-ab, \hspace{0.1in} \big (  A \big
  )_{i,j}=0 \hspace{0.1in} \text{otherwise}.
\end{equation*}
   The eigenvalues of $A$ are shown in \cite{Lo}  to be $a^2+b^2+2ab\cos(\frac {k\pi}{N+1}), 1 \leq k \leq N$, from which Theorem \ref{thm3.1} follows easily.

   We note here one can also give a proof of Theorem \ref{thm3.1} following the methods in  \cite{R2} and \cite{A1} as one checks readily that the right-hand side inequality of \eqref{3.7} follows on taking $\alpha=a, \beta=b, \mu_n=a^2+ab\sin(n+1)t/\sin(nt), t = \pi/(N+1)$ in inequality \eqref{2.2} and summing for $n=1, \ldots, N-1$.  Similarly,
  the left-hand side inequality of \eqref{3.7} follows on taking $\alpha=a, \beta=b, \mu_n=a^2-ab\sin(n+1)t/\sin(nt), t = \pi/(N+1)$ in inequality \eqref{2.2} (with inequality reversed there).

\section{Applications of Theorem \ref{thm03}}
\label{sec 4} \setcounter{equation}{0}

    As an application of Theorem \ref{thm03} or rather, Corollary \ref{cor03}, we now prove the following
\begin{theorem}
\label{thm04}
   Let $p \geq 2$ be fixed, then inequalities \eqref{8} hold for $1 \leq \alpha \leq 2$.
\end{theorem}
\begin{proof}
   For simplicity, we make a change of variable $\alpha-1 \mapsto \alpha$ so that by Corollary \ref{cor03}, it suffices to show the following inequality holds for any integer $n \geq 1$ and $0\leq \alpha \leq 1$ :
\begin{equation*}
  \frac {\sum^{n+1}_{k=1}k^{\alpha}}{(n+1)^{\alpha}}-\frac {\sum^{n}_{k=1}k^{\alpha}}{n^{\alpha}} \leq \frac 1{\alpha+1}+\frac {n^{\alpha}}{2\sum^{n}_{k=1}k^{\alpha}}\Big (1-\frac 1{p} \Big ) \frac 1{(\alpha+1)^2}.
\end{equation*}
  It is easy to see on letting $x_n=1/\sum^{n}_{k=1}k^{\alpha}$ that the above inequality is equivalent to $f_n(x_n) \geq 0$, where for $0 \leq x \leq 1$, we define
\begin{equation*}
  f_n(x)=1+n^{\alpha}x\Big (\frac 1{\alpha+1}+\frac {n^{\alpha}}{2}\Big(1-\frac 1{p}\Big) \frac x{(\alpha+1)^2} \Big )-\frac {n^{\alpha}}{(n+1)^{\alpha}}\Big ( 1+(n+1)^{\alpha}x \Big ).
\end{equation*}
   Now we need two lemmas:
 \begin{lemma}
 \label{lem10}    For $p \geq 2$ and $0 \leq \alpha \leq 1$, we have
 \begin{equation*}
    \frac 1{2^{\alpha}} \leq \frac 1{\alpha+1}+\frac {1-1/p}{2(\alpha+1)^2}.
 \end{equation*}
 \end{lemma}
 \begin{proof}
    As $p \geq 2$, we have $1-1/p \geq 1/2$. Hence our assertion is a consequence of the following inequality:
 \begin{equation*}
    \frac 1{2^{\alpha}} \leq \frac 1{\alpha+1}+\frac {1}{4(\alpha+1)^2}.
 \end{equation*}
    It is easy to see that the above inequality is equivalent to $h(\alpha) \geq 0$ for  $0 \leq \alpha \leq 1$,
    where
 \begin{equation*}
    h(\alpha)=2^{\alpha}(5+4\alpha)-4(1+\alpha)^2.
 \end{equation*}
    Note that by Taylor expansion, we have, for $0 \leq \alpha \leq 1$,
\begin{equation*}
   2^{\alpha}=(1+1)^{\alpha} \geq 1+\alpha +\alpha (\alpha-1)/2.
\end{equation*}
   It follows that
\begin{equation*}
   h(\alpha) \geq (1+\alpha +\alpha (\alpha-1)/2)(5+4\alpha)-4(1+\alpha)^2=(1-\alpha)(1-\alpha/2)+2\alpha^3 \geq 0.
\end{equation*}
     This now completes the proof.
 \end{proof}
   The above lemma implies that $f_1(x_1)=f_1(1) \geq 0$ for $p \geq 2$, $0 \leq \alpha \leq 1$ so now we may assume $n \geq 2$ and we need the following
\begin{lemma}\cite[Lemma 1, 2, p.18]{L&S}
\label{lem0}
    For an integer $n \geq 1$ and $0 \leq r \leq 1$,
\begin{equation}
\label{6.1}
    \frac {1}{r+1}n(n+1)^r \leq  \sum^n_{i=1}i^r \leq \frac {r}{r+1}\frac
   {n^r(n+1)^r}{(n+1)^r-n^r}.
\end{equation}
\end{lemma}
  The above lemma implies that $x_n \leq (\alpha+1)/(n(n+1)^{\alpha})$. Note that for fixed $n$, $f_n(x)$ is a quadratic function of $x$ and the only root of $f'_n(x)=0$ is $\alpha(\alpha+1)/(n^{\alpha}(1-1/p))$. Suppose we have
\begin{equation}
\label{7.2}
  \frac {\alpha+1}{n(n+1)^{\alpha}} \geq \frac {\alpha(\alpha+1)}{n^{\alpha}(1-1/p  )}.
\end{equation}
  Then it suffices to show that for fixed $0 \leq \alpha \leq 1$ and any $n$,
\begin{equation*}
  f_n \Big (\frac {\alpha(\alpha+1)}{n^{\alpha}(1-1/p)} \Big )= 1-\frac {\alpha^2}{2(1-1/p)}-\frac {n^{\alpha}}{(n+1)^{\alpha}} \geq  0.
\end{equation*}
  Note that \eqref{7.2} implies that
\begin{equation*}
  \frac {n^{\alpha}}{n(n+1)^{\alpha}} \geq \frac {\alpha}{1-1/p}.
\end{equation*}
  It follows that
\begin{equation*}
  f_n \Big (\frac {\alpha(\alpha+1)}{n^{\alpha}(1-1/p)} \Big ) \geq 1-\frac {\alpha n^{\alpha}}{2n(n+1)^{\alpha}} -\frac {n^{\alpha}}{(n+1)^{\alpha}}  .
\end{equation*}
  Thus it suffices to show the right-hand side expression is no less than $0$, which is equivalent to
\begin{equation}
\label{7.4}
  \Big (1+\frac {1}{n} \Big )^{\alpha}  \geq 1+\frac {\alpha}{2n}.
\end{equation}
   Note that by Taylor expansion, we have, for $0 \leq \alpha \leq 1$,
\begin{equation*}
   \Big (1+\frac {1}{n} \Big )^{\alpha} \geq 1+\frac {\alpha}{n} +\frac {\alpha (\alpha-1)}{2n^2}.
\end{equation*}
   Apply the above estimation in \eqref{7.4}, we see that inequality \eqref{7.4} will follow as long as $n \geq 1-\alpha$, which is certainly true.

   It remains to consider the case where inequality \eqref{7.2} reverses and we then deduce that when $n \geq 2$ and $0 \leq \alpha \leq 1$,
\begin{equation*}
  f_n(x_n) \geq f_n \Big (\frac {\alpha+1}{n(n+1)^{\alpha}} \Big )=\frac {n^{\alpha}}{(n+1)^{\alpha}}g(\frac 1{n}),
\end{equation*}
   where
\begin{equation*}
   g(y)=(1+y)^{\alpha}+y\Big (1+ \frac {1-1/p}{2}y(1+y)^{-\alpha} \Big )-1-(\alpha+1)y.
\end{equation*}
    Note that when $0 \leq \alpha  \leq 1$, $y>0$,
\begin{equation*}
   (1+y)^{\alpha} \geq 1+\alpha y+\alpha (\alpha-1)y^2/2; \hspace{0.1in} (1+y)^{-\alpha} \geq 1-\alpha y.
\end{equation*}
   We conclude from the above estimations that when $0 \leq y \leq 1/2$, $p \geq 2$, $0 \leq \alpha \leq 1$,
\begin{eqnarray*}
   g(y) &\geq & 1+\alpha y+\alpha (\alpha-1)y^2/2+y\Big (1+ \frac {1-1/p}{2}y(1-\alpha y) \Big )-1-(\alpha+1)y \\
   &=& \frac {y^2}{2}\Big (\alpha (\alpha-1)+(1-1/p)(1-\alpha y) \Big ) \geq \frac {y^2}{2}\Big (\alpha (\alpha-1)+(1-1/p)(1-\alpha/2) \Big ) \\
   & \geq & \frac {y^2}{2}\Big (\alpha (\alpha-1)+1/2(1-\alpha/2) \Big ) \geq 0.
\end{eqnarray*}
  This now implies that $f_n(x_n) \geq 0$ which in turn completes the proof.
\end{proof}

  We point out here that one can easily deduce from the proof of Theorem \ref{thm04} that inequalities \eqref{8} hold for $p \geq p_0$, $0 \leq \alpha \leq 1$ for a number $p_0$ with $1< p_0 <2$ or inequalities \eqref{8} hold for $p \geq p_1(\alpha)$, for a number $p_1(\alpha)$ depending on $\alpha, 0 \leq \alpha \leq 1$. But one can also see that the proof will not allow us to prove inequalities \eqref{8} for all $p>1$ and $0 \leq \alpha \leq 1$.

  Note that Theorem \ref{thm04} immediately implies, by our discussions in Section \ref{sec 1}, that inequality \eqref{1} holds for $\lambda_k=k^{\alpha}$ for $0 \leq \alpha \leq 1 $ with $E=e^{1/(\alpha+1)}$, which fills in the open cases of Bennett's conjecture. We now summarize this in the following
\begin{theorem}
\label{thm2}
   Inequality \eqref{1} holds for $\lambda_k=k^{\alpha}$ for $\alpha >-1 $ with $E=e^{1/(\alpha+1)}$.
\end{theorem}



\end{document}